\documentstyle[leqno,12pt,amstex]{article}
\numberwithin{equation}{section}
\begin{document}
\title
{
\large{\bf Regularity and lifespan of small solutions 
to systems of quasi-linear wave equations 
with multiple speeds, I: almost global existence}\\
{~}\\
}
\author
{
\large{Kunio Hidano}
}
\date{}
\maketitle
\begin{abstract}
In this paper, we show almost global existence 
of small solutions to the Cauchy problem 
for symmetric system of wave equations 
with quadratic (in $3$D) or cubic (in $2$D) 
nonlinear terms and multiple propagation speeds. 
To measure the size of initial data, 
we employ a weighted Sobolev norm whose 
regularity index is the smallest among all the admissible 
Sobolev norms of integer order. 
We must overcome the difficulty caused by the absence of 
the $H^1$-$L^p$ Klainerman-Sobolev type inequality, 
in order to obtain a required a priori bound in 
the low-order Sobolev norm. 
The introduction of good substitutes for this inequality 
is therefore at the core of this paper.  
Using the idea of showing 
the well-known Lady{\v z}enskaja inequality, 
we prove some weighted inequalities, 
which, together with the generalized Strauss inequality, 
play a role as the good substitute.
\end{abstract}


\section{Introduction}
\baselineskip=0.6cm
Let us start with some well-known results on the Cauchy problem for 
the quasi-linear, {\it scalar} wave equation in three space dimensions 
of the form
\begin{equation}
\partial_t^2 u-\Delta u
=
\sum_{\alpha,\beta,\gamma=0}^3
G^{\alpha\beta\gamma}(\partial_\alpha u)\partial^2_{\beta\gamma}u
+
\sum_{\alpha,\beta=0}^3
H^{\alpha\beta}(\partial_\alpha u)\partial_\beta u
,\,\,
t>0,\,x\in{\mathbb R}^3, 
\end{equation}
where $x_0:=t$, $\partial_\alpha:=\partial/\partial x_\alpha$, 
$\partial^2_{\alpha\beta}:=\partial^2/\partial x_\alpha\partial x_\beta$, 
and 
$G^{\alpha\beta\gamma}$, $H^{\alpha\beta}$ 
are real constants. 
Though our primary concern in the present paper is on 
the Cauchy problem for the {\it system} of nonlinear wave equations 
with multiple speeds, 
we expect that to revisit some fundamental results for 
the scalar equation (1.1) 
will serve as a guide to the main problem 
discussed later. 
Just for simplicity, we suppose that 
$G^{\alpha00}=0$ for any $\alpha$. 
Moreover, without loss of generality, 
we may suppose $G^{\alpha\beta\gamma}=G^{\alpha\gamma\beta}$ 
for any $\alpha, \beta$, and $\gamma$ 
because we are interested in classical solutions. 
It is then well known 
(see, e.g., page 113 of H\"ormander \cite{Hor}) 
that for any initial data in $H^4\times H^3$ 
with 
$$
\sum_{\alpha=0}^3
\sum_{b,c=1}^3
|G^{\alpha bc}|\|\partial_\alpha u(0)\|_{L^\infty}
$$ 
small enough, 
there exists $T>0$ and 
a unique solution 
$u$ to the Cauchy problem for (1.1) such that 
$\partial_\alpha u
\in
C([0,T];H^3)$. 
If we consider (1.1) in ${\mathbb R}^n$ 
instead of ${\mathbb R}^3$ 
and work with fractional-order Sobolev spaces, 
this result of local existence remains true for 
the $H^{s+1}\times H^s$-data with $s>(n/2)+1$ 
(see, e.g, Proposition 5.2.B of Taylor \cite{Tay}). 
Note that, for $n=2,3$, 
$H^4\times H^3$ is the largest 
among all the admissible Sobolev spaces of integer order, 
as far as the standard local existence theorem is concerned. 
We also remark that, 
when considering (1.1) in ${\mathbb R}^n$ ($n=2,3$) 
with more regular data in $H^{s+1}\times H^s$ 
for some $s\geq 4$ $(s\in{\mathbb N})$, 
it is possible to choose $T>0$ depending only on 
$\sum_\alpha\|\partial_\alpha u(0)\|_{H^3}$ 
and independently of $s$ 
such that the Cauchy problem for (1.1) admits 
a unique solution satisfying 
$\partial_\alpha u\in C([0,T];H^s)$. 
(See, e.g., Theorem 5.8 of Racke \cite{Rac}. 
Note that the equation (1.1) can be written in the form 
of the first order quasi-linear system (see, e.g., page 19 
of John \cite{JohnULS}) to which we can apply 
Proposition 5.2.B of \cite{Tay} 
and Theorem 5.8 of \cite{Rac}.) 
This means that 
a continuation of local smooth solutions to a larger strip 
is reduced to the a priori $H^3$-bound of 
their first derivatives. 

Concerning long-time existence, 
there exist 
positive constants $C_1$, $\varepsilon_1$ depending on 
the coefficients $G^{\alpha\beta\gamma}$, $H^{\alpha\beta}$ 
such that 
whenever compactly supported $C^\infty$-initial data satisfies 
\begin{equation}
{\cal W}_4^{1/2}(u(0))\leq\varepsilon_1, 
\end{equation}
we can obtain the a priori $L^2$-bound 
of $\partial_\alpha\Gamma^a u(t)$, 
$\alpha=0,\dots,3,\, |a|\leq 3$ 
which is strong enough to show that 
the smooth local solution exists at least for the interval 
$[0,\exp(C_1\varepsilon^{-1})]$, 
where $\varepsilon:={\cal W}_4^{1/2}(u(0))$. 
Here and in the following discussion, we use the notation: 
\begin{eqnarray}
& &
{\cal W}_1(u(t))
:=
\frac12
\int_{{\mathbb R}^3}
\bigl(
(\partial_t u(t,x))^2
+
|\nabla u(t,x)|^2
\bigr)dx,\\
& &
{\cal W}_\kappa(u(t))
:=
\sum_{|a|\leq \kappa-1}
{\cal W}_1(\Gamma^a u(t)),
\quad \kappa=2,3,\dots
\end{eqnarray}
For a multi-index $a$, 
$\Gamma^a$ stands for any product 
of the $|a|$ operators 
$\partial_\alpha$ $(\alpha=0,\dots,3)$, 
$\Omega_{ij}:=x_i\partial_j-x_j\partial_i$ 
$(1\leq i<j\leq 3)$, 
$L_k:=x_k\partial_t+t\partial_k$ 
$(k=1,2,3)$, 
and $S:=t\partial_t+x\cdot\nabla$. 
The proof of the a priori bound 
uses the Sobolev-type inequality 
\begin{equation}
\|v(t,\cdot)\|_{L^p({\mathbb R}^3)}
\leq
C(1+t)^{-2(1/2-1/p)}
\sum_{|a|\leq 1}
\|\Gamma^a v(t,\cdot)\|_{L^2({\mathbb R}^3)}
\end{equation}
(see, e.g., Ginibre and Velo \cite{GV}) 
as well as the Klainerman inequality \cite{Kl1987}
\begin{equation}
\|v(t,\cdot)\|_{L^\infty({\mathbb R}^3)}
\leq
C(1+t)^{-1}
\sum_{|a|\leq 2}
\|\Gamma^a v(t,\cdot)\|_{L^2({\mathbb R}^3)}.
\end{equation}
(Note that, while a loss of just one derivative occurs in (1.5), 
we lose two derivatives in applying (1.6) to the estimation 
of nonlinear terms.) 
We note that the interval mentioned above 
becomes exponentially large 
as the size $\varepsilon$ of initial data 
gets smaller and smaller. 
Such results have been called 
``almost global existence theorem'' 
in the literature 
since the pioneering work of 
John and Klainerman \cite{JK} for the equation (1.1). 

Now, let us turn our attention to 
the main concern in the present paper: 
the Cauchy problem for the system of 
nonlinear wave equations of the form
\begin{align}
(\partial_t^2-c_l^2\Delta)u^l
&=
G_{ij}^{l,\alpha\beta\gamma}
(\partial_\alpha u^i)\partial^2_{\beta\gamma}u^j\\
&+
H_{ij}^{l,\alpha\beta}(\partial_\alpha u^i)\partial_\beta u^j
,\,\,
t>0,\,x\in{\mathbb R}^3\nonumber
\end{align}
and its $2$D counterpart
\begin{align}
(\partial_t^2-c_l^2\Delta)u^l
&=
G_{ijk}^{l,\alpha\beta\gamma\delta}
(\partial_\alpha u^i)(\partial_\beta u^j)
\partial_{\gamma\delta}^2 u^k\\
&+
H_{ijk}^{l,\alpha\beta\gamma}(\partial_\alpha u^i)(\partial_\beta u^j)
\partial_\gamma u^k
,\,\,
t>0,\,x\in{\mathbb R}^2.\nonumber
\end{align}
Here, 
$u=(u^1,\dots,u^N):(0,T)\times{\mathbb R}^n
\to{\mathbb R}^N$ 
$(n=2,3,\,N\in{\mathbb N})$ 
and, on the right-hand side of (1.7)--(1.8) 
and in the following discussion as well, 
repeated indices are summed 
if lowered and uppered. 
Greek indices range from $0$ to $n$, 
and roman indices from $1$ to $N$. 
Suppose the symmetry condition
\begin{equation}
G_{ij}^{l,\alpha\beta\gamma}
=
G_{ij}^{l,\alpha\gamma\beta}
=
G_{il}^{j,\alpha\beta\gamma}
\end{equation}
for any $i,j,l$ and $\alpha, \beta, \gamma$ 
in (1.7) or 
\begin{equation}
G^{l,\alpha\beta\gamma\delta}_{ijk}
=
G^{l,\alpha\beta\delta\gamma}_{ijk}
=
G^{k,\alpha\beta\gamma\delta}_{ijl}
\end{equation} 
for any $i,j,k,l$ and $\alpha, \beta, \gamma, \delta$ 
in (1.8). 
Then the results of local existence mentioned above for (1.1) 
carry over to (1.7), (1.8), 
because these systems can be written in the form 
of the first order quasi-linear symmetric system, 
such as (5.9) of \cite{Rac}, (5.2.1) of \cite{Tay}. 
Concerning long-time existence, 
for compactly supported smooth initial data of the form
\begin{equation}
u(0)=\varepsilon f,\,\,\partial_t u(0)=\varepsilon g,
\end{equation}
Sogge proved that 
there exist positive constants $C_2$, $\varepsilon_2$ 
depending on the speeds $c_1,\dots,c_N$, 
the coefficients $G_{ij}^{l,\alpha\beta\gamma}$, 
$H_{ij}^{l,\alpha\beta}$, and 
a weighted $H^{10}$-norm of $\nabla f$ and $g$ such that 
if $0<\varepsilon\leq\varepsilon_2$, 
then a unique solution to (1.7), (1.11) exists 
at least over the time interval 
$[0,\exp(C_2\varepsilon^{-1})]$. 
(See Theorem 4.1 on page 67 of Sogge \cite{So2ndEd}, 
the proof of which is based on that of 
Theorem 1.2 of Keel, Smith, and Sogge \cite{KSS-JAMS}.) 
Also, Kovalyov \cite{Kov} proved that 
there exist positive constants $C_3$, $\varepsilon_3>0$ 
depending on the speeds $c_1,\dots,c_N$, 
the coefficients $G_{ijk}^{l,\alpha\beta\gamma\delta}$, 
$H_{ijk}^{l,\alpha\beta\gamma}$, and 
a weighted $H^5$-norms of 
$\nabla f$ and $g$ such that 
if $0<\varepsilon\leq\varepsilon_3$, 
then a unique solution to (1.8), (1.11) 
exists at least over 
$[0,\exp(C_3\varepsilon^{-2})]$. 
In the present paper, 
we aim at refining these results 
by employing a lower-order norm 
to measure the size of initial data. 
More precisely, we prove:

\smallskip

\noindent{\bf Theorem 1.1.} {\it Assume 
the symmetry condition $(1.9)$, $(1.10)$. 
Then there exist positive constants 
$C_0$, $\varepsilon_0$ depending on the propagation speeds and 
the coefficients of the equations $(1.7)$, $(1.8)$ such that 
if compactly supported, smooth initial data is small so that
\begin{equation}
N_4(u(0))\leq\varepsilon_0
\end{equation}
may hold, then the systems $(1.7)$, $(1.8)$ admit unique solutions 
defined on the interval $[0,T]$ such that 
$N_4(u(t))\leq 2N_4(u(0))$, $0\leq t\leq T$. 
Here 
$$
T=\exp(C_0\varepsilon^{-\nu})
\quad
(\varepsilon:=N_4(u(0)))
$$
with $\nu=1$ for $(1.7)$, 
$\nu=2$ for $(1.8)$.}

\smallskip

Here, on the basis of the standard energy $E_1(u(t))$ 
associated with unperturbed wave equations 
\begin{equation}
E_1(u(t))
=
\frac12
\sum_{l=1}^N
\int_{{\mathbb R}^n}
\bigl(
|\partial_t u^l(t,x)|^2
+
c_l^2|\nabla u^l(t,x)|^2
\bigr)
dx,
\end{equation}
we have defined the quantity $N_\kappa(v(t))$ 
for $v=(v^1,\dots,v^N)$ as 
\begin{align}
&
N_1(v(t))
=
\sqrt{E_1(v(t))},\,
N_2(v(t))
=
\left(
\sum_{|a|+|b|+d\leq 1}
E_1(\partial_x^a \Omega^b S^d v(t))
\right)^{1/2},\\
&
N_\kappa(v(t))
=
\left(
\sum_{{|a|+|b|+d\leq \kappa-1}\atop{d\leq 1}}
E_1(\partial_x^a \Omega^b S^d v(t))
\right)^{1/2}
,\,\kappa=3,4,\dots,\nonumber
\end{align}
where, for $a=(a_1,\dots,a_n)$ and $b=(b_1,\dots,b_m)$ 
($m=1,3$ for $n=2,3$, respectively), 
$\partial_x^a \Omega^b S^d v
:=
(\partial_x^a \Omega^b S^d v^1,\dots,\partial_x^a \Omega^b S^d v^N)$, 
$\partial_x^a:=\partial_1^{a_1}\cdots\partial_n^{a_n}$, 
$\Omega^b:=\Omega_{12}^{b_1}\cdots\Omega_{n-1\,n}^{b_m}$, 
$\Omega_{ij}=x_i\partial_j-x_j\partial_i$ 
$(1\leq i<j\leq n)$, 
and $S=t\partial_t+x\cdot\nabla$. 
We set 
$Z:=\{\partial_i,\Omega_{jk},S:i=1,\dots,n, 1\leq j<k\leq n\}$. 
Note that none of the operators 
$L_k=x_k\partial_t+t\partial_k$ $(k=1,\dots,n)$ 
is an element of the set $Z$. 
Note also $\partial_t\notin Z$. 

There exist some difficulties in showing 
the almost global existence result for (1.7), (1.8) 
when we employ the lower-order norm, such as $N_4$ defined above, 
to measure the size of data. 
Recall that, besides the standard energy inequality 
for the variable-coefficient wave equation, 
the generalized Sobolev-type inequalities (1.5)--(1.6) 
and the nice commutation relations between 
the D'Alembertian (with the propagation speed $c=1$) 
and the elements of $\{\Omega_{ij},L_k,S\}$ 
play an important role 
in showing the almost global existence for (1.1) 
with compactly supported, smooth data satisfying (1.2). 
When considering the multiple-speed system 
(1.7), (1.8), 
we must take into account the fact that 
the operator ${\tilde L}_k:=c^{-1}x_k\partial_t+ct\partial_k$, 
which is a speed-dependent variant of $L_k$, 
commutes with 
$\Box_c:=\partial_t^2-c^2\Delta$, 
but it no longer does with $\Box_{\hat c}:=\partial_t^2-{\hat c}^2\Delta$ 
$({\hat c}\ne c)$. 
Indeed, we have 
$[{\tilde L}_k,\Box_{\hat c}]
=2c^{-1}({\hat c}^2-c^2)\partial_k\partial_t$, 
and this commutation relation is obviously useless 
in our argument. 
We must therefore give up using such a modified operator 
${\tilde L}_k$, 
which in turn means that we must give up using the Sobolev-type 
inequalities (1.5)--(1.6). 
On the other hand, 
we still enjoy the good commutation relations 
$[\Omega_{ij},\Box_c]=0$ 
and $[S,\Box_c]=-2\Box_c$, 
and some good substitutes for the Klainerman inequality (1.6) 
are available on the basis of the use of 
the operators $\Omega_{ij}$ and $S$ 
and without relying upon the operators ${\tilde L}_k$. 
(See Lemma 6.1 of \cite{ST} and Lemma 1 of \cite{Tom2}. 
See also (4.2) of \cite{HiTohoku}.) 
These substitutes, combined with 
the Klainerman-Sideris inequality (see (3.1) below), 
would suffice to show almost global existence theorem 
for (1.7), (1.8) when a suitable higher (than 4) order norm of data is small enough. 
See, e.g., Section 8 of Sideris and Tu \cite{ST} and 
Theorem 3.1 of \cite{HiTohoku}. 
Therefore, it is a good substitute for (1.5) 
that plays a key role in 
reducing the regularity index of norm to as low a level as in (1.2). 
To the best of the present author's knowledge, 
no substitute for (1.5) is available in the literature. 
We explain that our key weighted inequalities (2.7)--(2.13) 
are well combined with the method of \cite{ST}, 
and they play a role as the substitute for (1.5). 
To prove these key inequalities, 
we follow the way of showing 
the well-known Lady{\v z}enskaja inequality \cite{Lad} 
or use the generalization of the Strauss inequality. 
(See (2.20) for the generalized Strauss inequality. 
See also \cite{CO} for recent, another extension 
of the classical inequality of Strauss \cite{Strauss}.) 

The method in this paper has an application 
to the system of quasi-linear wave equations with 
{\it quadratic} nonlinear terms in $2$D. 
Repeating essentially the same argument as in the proof of 
Theorem 1.1, we obtain the following:

\smallskip

\noindent{\bf Theorem 1.2.} {\it Consider $(1.7)$ in ${\mathbb R}^2$. 
Assume 
the symmetry condition $(1.9)$. 
Then there exist positive constants 
$A_0$, $\varepsilon_0$ depending on the propagation speeds and 
the coefficients of $(1.7)$ 
with the following property$:$ 
if the compactly supported, smooth initial data is small so that 
$N_4(u(0))\leq\varepsilon_0$ may hold, 
then the problem $(1.7)$ has a unique solution satisfying 
$N_4(u(t))\leq 2N_4(u(0))$,\,\,$0<t<T$. 
Here 
$T=A_0\varepsilon^{-2}$ 
$(\varepsilon:=N_4(u(0)))$.
}

\smallskip

In \cite{Kov}, Kovalyov considered the system (1.7) 
not in ${\mathbb R}^3$ but in ${\mathbb R}^2$ with data of the form (1.11), 
and obtained the slightly weaker lower bound 
$T\geq C\varepsilon^{-2}(\log (1/\varepsilon))^{-2}$, 
while in \cite{HosAMSA}, 
assuming $H_{ij}^{l,\alpha\beta}=0$ for all $i,j,l,\alpha,\beta$, 
Hoshiga obtained the refined lower bound 
$T\geq C\varepsilon^{-2}$ 
with a positive constant $C$ 
computed explicitly from 
the propagation speeds $c_1,\dots,c_N$, 
the coefficients $G_{ij}^{l,\alpha\beta\gamma}$, 
and the given functions $f$, $g$. 
Theorem 1.2 is an improvement on the previous results 
of \cite{Kov} and \cite{HosAMSA}, for 
Kovalyov used a higher-order norm to 
measure the size of initial data and 
his lower bound of the lifespan is slightly weaker than ours, 
and Hoshiga imposed the restriction 
$H_{ij}^{l,\alpha\beta}=0$ for all $i,j,l,\alpha,\beta$, 
while we no longer need his strict restriction. 

Here we give three remarks. 
Firstly, as in the books \cite{Al}, \cite{Hor}, \cite{JohnULS}, 
\cite{Rac}, and \cite{So2ndEd}, 
we have so far supposed that initial data is smooth and 
compactly supported, 
when considering the lifespan of small solutions. 
This is mainly because a continuation argument becomes 
considerably easier for compactly supported 
(in space at fixed times $t>0$), smooth solutions. 
See (4.31) below. 
Note that the constants $A_0$, $C_0$, and $\varepsilon_0$ 
appearing Theorems 1.1 and 1.2 are completely independent of 
the ``size'' of the support of initial data. 
Therefore, once we have proved these theorems, 
we should move on to removing the compactness assumption of 
the support, as well as the regularity ($C^\infty$) 
assumption, of initial data. 
The idea of doing it can be found on page 122 
of \cite{Hor} (see {\it Remark} there). 
In order to keep the present paper to a moderate length, 
we refrain from pursuing this important problem.

Secondly, in the definition of $N_\kappa(u(t))$ 
$(\kappa\geq 3)$ 
we have limited the number of occurrences of $S$ to $1$, 
in accordance with the idea of 
the earlier papers \cite{KSS-JAMS} and \cite{HY2} that 
its at most $1$ occurrence is 
actually sufficient for the proof of almost global existence. 
With this, there is an advantage that 
we can bypass the burdensome calculation of $\partial_t^j u(0,x)$ $(j=2,3,4)$ 
when computing $N_4(u(0))$, 
because $\partial_t\notin Z$ 
and $\partial_tSu=\partial_tu+x\cdot\nabla\partial_t u$ 
at $t=0$. 
(Compare this with the fact that 
we must successively calculate $\partial_t^j u(0,x)$ $(j=2,3,4)$ 
with the help of the equation (1.1) 
when computing ${\cal W}_4^{1/2}(u(0))$ appearing in (1.2).) 
Another feature lies in that, 
when initial data 
$(u(0),\partial_t u(0))=(\varphi,\psi)$ 
is radially symmetric about $x=0$ 
(and the system of equations is not necessarily so), 
we easily see the condition (1.12) 
is satisfied whenever 
the norm with the ``mild'' weight 
$\langle x\rangle:=\sqrt{1+|x|^2}$ 
\begin{equation}
\sum_{1\leq l\leq N}
\bigl(
\sum_{1\leq |a|\leq 4}
\|\langle x\rangle\partial_x^a \varphi^l\|_{L^2}
+
\sum_{|a|\leq 3}
\|\langle x\rangle\partial_x^a \psi^l\|_{L^2}
\bigr)
\end{equation}
is small enough. 
The result of almost global existence 
for symmetric (and not necessarily diagonal) 
systems of quasi-linear wave equations is new 
when smallness is required of only such mildly weighted 
Sobolev norm of radial data. 

Thirdly, the proof of Theorem 1.1 obviously remains valid 
for the scalar equation (1.1), 
thus we obtain almost global existence result 
under the condition (1.12) with $N=1$ which is weaker than (1.2). 

We conclude this section by mentioning that, 
in the sequel \cite{sequel}, 
assuming the null condition 
in the different-speed setting proposed by 
Agemi and Yokoyama \cite{AY}, Yokoyama \cite{Yo}, 
we will prove the global existence theorem for (1.7), (1.8) 
on a condition which is stronger than (1.12), 
but weaker than that in the previous papers \cite{Yo}, \cite{HK}, 
\cite{ST}, \cite{So2ndEd}, \cite{LNS}.
In addition to the key tools used in the present paper, 
the proof will use the estimation lemmas due to 
Sideris and Tu (see Lemma 5.1 of \cite{ST}) in $3$D, 
Lindblad, Nakamura, and Sogge (see Lemma A.4 of \cite{LNS}) 
in $2$D when handling the null-form terms. 

This paper is organized as follows. 
In the next section, some useful inequalities of 
the Sobolev type or the trace type are proved. 
Using the Klainerman-Sideris inequality, 
we bound weighted space-time $L^2$-norms of 
the second or some higher-order derivatives 
of the local solution in Section 3. 
In Sections 4 and 5, we carry out 
the energy integral argument 
and complete the proof of Theorems 1.1 and 1.2. 
\section{Preliminaries}
As explained in Section 1, 
repeated indices will be summed if lowered and uppered. 
Greek indices range from $0$ to $n$ ($n=2$ or $3$), 
and roman indices from $1$ to $N$ or $1$ to $n$. 
In addition to the usual partial differential operators 
$\partial_\alpha=\partial/\partial x_\alpha$ 
$(\alpha=0,\dots,n)$, 
we use the generator of Euclid rotation 
$\Omega_{ij}=x_i\partial_j-x_j\partial_i$ 
and of space-time scaling $S=t\partial_0+x\cdot\nabla$. 
The set of these $\mu$ 
($\mu=4$ for $n=2$, $\mu=7$ if $n=3$) 
differential operators is denoted by 
$Z=\{\,Z_1,\dots,Z_\mu\,\}
=\{\,\nabla,\Omega,S\,\}$. 
Note that $\partial_t$ is not an element of $Z$. 
For a multi-index $a=(a_1,\dots,a_\mu)$, 
we set $Z^a:=Z_1^{a_1}\cdots Z_\mu^{a_\mu}$. 
We also use  
${\bar Z}=\{\,Z_1,\dots,Z_{\mu-1}\,\}
=\{\,\nabla,\Omega\,\}$, 
with 
${\bar Z}^a:=\partial_1^{a_1}\partial_2^{a_2}\Omega_{12}^{a_3}$ 
$(a=(a_1,a_2,a_3))$, 
${\bar Z}^a:=\partial_1^{a_1}\partial_2^{a_2}\partial_3^{a_3}
\Omega_{12}^{a_4}\Omega_{13}^{a_5}\Omega_{23}^{a_6}$ 
$(a=(a_1,\dots,a_6))$ 
for $n=2,3$, respectively. 

We collect several results concerning 
commutation relations 
and Sobolev-type and trace-type inequalities. 
Let $[\cdot,\cdot]$ be the commutator: 
$[A,B]:=AB-BA$. 
It is easy to verify that 
\begin{eqnarray}
& &
[Z_i,\partial_t^2-c^2\Delta]=0\,\,\,\mbox{for $i=1,\dots,\mu-1$},\,\,\,
[S,\partial_t^2-c^2\Delta]=-2(\partial_t^2-c^2\Delta),\\
& &
[Z_j,Z_k]=\sum_{i=1}^\mu C^{j,k}_i Z_i,\,\,\,
j,\,k=1,\dots,\mu,\\
& &
[Z_j,\partial_k]
=
\sum_{i=1}^n C^{j,k}_i\partial_i,\,\,\,j=1,\dots,\mu,\,\,k=1,\dots,n,\\
& &
[Z_j,\partial_t]=0,\,j=1,\dots,\mu-1,\quad [S,\partial_t]=-\partial_t.
\end{eqnarray}
Here $C^{j,k}_i$ denotes a constant depending on 
$i$, $j$, and $k$.

The following lemma is concerned with Sobolev-type 
or trace-type inequalities. 
We use these inequalities in combination with 
the Klainerman-Sideris inequality (see (3.1) below). 
The auxiliary norms of $v=(v^1,\dots,v^N)$
\begin{align}
M_2(v(t))
&=\sum_{l=1}^N
\sum_{{0\leq\delta\leq n}\atop{1\leq j\leq n}}
\|\langle c_l t-|x|\rangle
\partial_{\delta j}^2 v^l(t)\|_{L^2({\Bbb R}^n)},\\
M_4(v(t))
&=\sum_{|a|\leq 2}M_2({\bar Z}^a v(t)),
\nonumber
\end{align}
which appear in the following discussion, play an intermediate role. 
We remark that $S$ and $\partial_t^2$ are absent in the 
right-hand side above. 
Here and later on as well, 
we use the standard notation 
$\langle A\rangle=\sqrt{1+|A|^2}$ for a scalar or a vector $A$. 
We also use the notation $\partial_r:=(x/|x|)\cdot\nabla$, 
\begin{align}
&
\|w\|_{L_r^\infty L_\omega^p({\mathbb R}^n)}
:=
\sup_{r>0}
\|w(r\cdot)\|_{L^p(S^{n-1})},\\
&
\|w\|_{L_r^2 L_\omega^p({\mathbb R}^n)}
:=
\biggl(
\int_0^\infty \|w(r\cdot)\|_{L^p(S^{n-1})}^2 r^{n-1}dr
\biggr)^{1/2}.\nonumber
\end{align}

\smallskip

\noindent{\bf Lemma 2.1.} {\it Let 
$v$ be a vector-valued function 
$v=(v^1,\dots,v^N):(0,\infty)\times{\mathbb R}^n\to{\Bbb R}^N$ 
decaying sufficiently fast as $|x|\to\infty$. 
The following inequalities hold for every $l=1,\dots,N:$ 

\noindent{\rm (i)} Suppose $n=2$. We have for $\alpha=0,1,2$
\begin{align}
&
\|r^{1/2}
\partial_\alpha v^l(t)\|_{L_r^\infty L_\omega^2({\mathbb R}^2)}
\leq
C
N_1^{1/2}(v(t))
\bigl(
\sum_{|a|=1}
N_1(\partial_x^a v(t))
\bigr)^{1/2},\\
&
\|\langle c_l t-r\rangle^{1/2}\partial_\alpha v^l(t)\|_{L^4({\mathbb R}^2)}
\leq
CN_1^{1/2}(v(t))
\bigl(
N_1(v(t))
+
M_2(v(t))
\bigr)^{1/2}.
\end{align}
\noindent{\rm (ii)} Suppose $n=3$. We have for $\alpha=0,1,2,3$
\begin{align}
&
\|\langle c_l t-r\rangle^{1/2}\partial_\alpha v^l(t)\|_{L^3({\mathbb R}^3)}
\leq
CN_1^{1/2}(v(t))
\bigl(
N_1(v(t))
+
M_2(v(t))
\bigr)^{1/2},\\
&
\|
\langle c_l t-r\rangle\partial_\alpha v^l(t)
\|_{L^6({\mathbb R}^3)}
\leq
C
\bigl(
N_1(v(t))
+
M_2(v(t))
\bigr).
\end{align}
Moreover, for any $2\leq p<4$ 
there exists a constant $C=C_p>0$ 
such that we have 
\begin{equation}
\|r\partial_\alpha v^l(t)\|_{L_r^\infty L_\omega^p({\mathbb R}^3)}
\leq
C
\sum_{|a|\leq 1}
N_1({\bar Z}^a v(t)).
\end{equation}
}

\smallskip

{\it Remark}. We give three remarks. 
Firstly, by the Sobolev embedding 
$W^{1,4}({\mathbb R}^2)\hookrightarrow L^\infty({\mathbb R}^2)$ 
and 
$W^{1,6}({\mathbb R}^3)\hookrightarrow L^\infty({\mathbb R}^3)$, 
we get from (2.8) and (2.10)
\begin{align}
&\langle c_l t-r\rangle^{1/2}|\partial_\alpha v^l(t,x)|\\
&
\leq
C
\biggl(
\sum_{|a|\leq 1}N_1(\partial_x^a v(t))
\biggr)^{1/2}
\biggl(
\sum_{|a|\leq 1}N_1(\partial_x^a v(t))
+
\sum_{|a|\leq 1}M_2(\partial_x^a v(t))
\biggr)^{1/2}\nonumber
\end{align}
for $n=2$ and 
\begin{equation}
\langle c_l t-r\rangle
|\partial_\alpha v^l(t,x)|
\leq
C
\biggl(
\sum_{|a|\leq 1}N_1(\partial_x^a v(t))
+
\sum_{|a|\leq 1}M_2(\partial_x^a v(t))
\biggr)
\end{equation}
for $n=3$, respectively. 
The former was shown by Sideris 
(see the last inequality on page 379 of \cite{Tom2}). 
After he submitted the manuscript, 
the author became aware of the recent paper of Zha \cite{Z} 
where the latter (2.13) had been proved (see (37) there). 
In addition to (2.8)--(2.10), 
we will also use both (2.12) and (2.13) in the following discussion. 
We also note that 
the multiplicative form of the right-hand side of (2.8) and (2.12) 
is very useful in our argument 
(see, e.g., (3.15)--(3.16) below). 
Secondly, we remark that we will also use for $n=2,3$
\begin{equation}
\|r^{(n-1)/2}\partial_\alpha v^l(t)\|_{L^\infty({\mathbb R}^n)}
\leq
C\sum_{|a|\leq 2}
N_1({\bar Z}^a v(t)),
\end{equation}
which follows immediately from 
the combination of (2.7), (2.11) 
with the Sobolev embedding 
$W^{1,2}(S^1)\hookrightarrow L^\infty(S^1)$, 
$W^{1,p}(S^2)\hookrightarrow L^\infty(S^2)$ with 
$p>2$, respectively. 
Thirdly, in fact we will use (2.9) not in the present paper 
but in \cite{sequel}. 
The proof of (2.9) is similar to that of (2.8), 
thus we prove it here.

\smallskip

{\it Proof of Lemma $2.1$}. Applying the well-known inequality 
$\|\varphi\|_{L^6({\mathbb R}^3)}\leq C\|\nabla\varphi\|_{L^2({\mathbb R}^3)}$ 
with $\varphi=\langle c_l t-r\rangle\partial_\alpha v^l(t,x)$, 
we easily obtain (2.10). 
The proof of (2.8) builds upon how to obtain the well-known 
Lady{\v z}enskaya inequality \cite{Lad} 
$\|\varphi\|_{L^4({\mathbb R}^2)}^4
\leq
4
\|\varphi\|_{L^2({\mathbb R}^2)}^2
\|\nabla\varphi\|_{L^2({\mathbb R}^2)}^2
$. 
Indeed, we first obtain by a direct computation
\begin{align}
&\langle c_lt-r\rangle
|\partial_\alpha v^l(t,x)|^2
=
\int_{-\infty}^{x_1}
\frac{d}{d\xi_1}
\bigl(
\langle c_lt-{\tilde r}\rangle|\partial_\alpha v^l(t,\xi_1,x_2)|^2
\bigr)d\xi_1\\
&
\leq
C\int_{-\infty}^\infty
\bigl(
|\partial_\alpha v^l(t,\xi_1,x_2)|^2
+
\langle c_lt-{\tilde r}\rangle
|\partial_\alpha v^l(t,\xi_1,x_2)|
|\partial^2_{1\alpha} v^l(t,\xi_1,x_2)|
\bigr)d\xi_1,\nonumber
\end{align}
which yields
\begin{align}
&\langle c_lt-r\rangle^2
|\partial_\alpha v^l(t,x)|^4\\
&
\leq
C\int_{-\infty}^\infty
\bigl(
|\partial_\alpha v^l(t,\xi_1,x_2)|^2
+
\langle c_lt-{\tilde r}\rangle
|\partial_\alpha v^l(t,\xi_1,x_2)|
|\partial^2_{1\alpha} v^l(t,\xi_1,x_2)|
\bigr)d\xi_1\nonumber\\
&
\times
\int_{-\infty}^\infty
\bigl(
|\partial_\alpha v^l(t,x_1,\xi_2)|^2
+
\langle c_lt-{\hat r}\rangle
|\partial_\alpha v^l(t,x_1,\xi_2)|
|\partial^2_{2\alpha} v^l(t,x_1,\xi_2)|
\bigr)d\xi_2,\nonumber
\end{align}
where ${\tilde r}:=(\xi_1^2+x_2^2)^{1/2}$, 
${\hat r}:=(x_1^2+\xi_2^2)^{1/2}$. 
Integrating both the sides above over ${\mathbb R}^2$ and 
using the Fubini theorem and the Schwarz inequality, 
we get
\begin{align}
&
\int_{{\mathbb R}^2}
\langle c_lt-r\rangle^2
|\partial_\alpha v^l(t,x)|^4dx\\
&
\leq
C
\biggl(
N_1(v(t))
\bigl(
N_1(v(t))+M_2(v(t))
\bigr)
\biggr)^2,\nonumber
\end{align}
as desired. 
The proof of (2.9) is similar, 
and we follow the proof of 
$\|\varphi\|_{L^3({\mathbb R}^3)}
\leq\sqrt{2}
\|\varphi\|_{L^2({\mathbb R}^3)}^{1/2}
\|\nabla\varphi\|_{L^2({\mathbb R}^3)}^{1/2}$ 
which is a special case of 
the Gagliardo-Nirenberg inequality 
(see, e.g., page 25 of \cite{Fried}). 
As in (2.16), we get
\begin{align}
&(\langle c_lt-r\rangle^{1/2}
|\partial_\alpha v^l(t,x)|)^3\\
&
\leq
C
\biggl(
\int_{-\infty}^\infty
\bigl(
|\partial_\alpha v^l(t,X_1)|^2
+
\langle c_lt-{\tilde r}\rangle
|\partial_\alpha v^l(t,X_1)|
|\partial^2_{1\alpha} v^l(t,X_1)|
\bigr)d\xi_1
\biggr)^{1/2}
\nonumber\\
&
\times
\biggl(
\int_{-\infty}^\infty
\bigl(
|\partial_\alpha v^l(t,X_2)|^2
+
\langle c_lt-{\hat r}\rangle
|\partial_\alpha v^l(t,X_2)|
|\partial^2_{2\alpha} v^l(t,X_2)|
\bigr)d\xi_2
\biggr)^{1/2}
\nonumber\\
&
\times
\biggl(
\int_{-\infty}^\infty
\bigl(
|\partial_\alpha v^l(t,X_3)|^2
+
\langle c_lt-{\bar r}\rangle
|\partial_\alpha v^l(t,X_3)|
|\partial^2_{3\alpha} v^l(t,X_3)|
\bigr)d\xi_3
\biggr)^{1/2},\nonumber
\end{align}
where $X_1:=(\xi_1,x_2,x_3),\dots,
X_3:=(x_1,x_2,\xi_3)$, 
${\tilde r}:=|X_1|$, 
${\hat r}:=|X_2|$, and 
${\bar r}:=|X_3|$. 
Integrating both the sides above over ${\mathbb R}^3$ and 
using the Schwarz inequality repeatedly, we get (2.9).

The other inequalities (2.7), (2.11) 
follow from the well-known inequality 
(called the Strauss inequality, especially when 
we focus on radially symmetric functions; see \cite{Strauss} and \cite{CO})
\begin{equation}
\|r^{(n-1)/2}\varphi\|_{L_r^\infty L_\omega^2({\mathbb R}^n)}
\leq
\sqrt{2}
\|\partial_r \varphi\|_{L^2({\mathbb R}^n)}^{1/2}
\|\varphi\|_{L^2({\mathbb R}^n)}^{1/2}
\end{equation}
or its generalization (see (2.10) of $\cite{HWYMathAnn}$): 
for $2\leq q\leq \infty$ and $2/p=1/2+1/q$ 
(the reader is asked 
to interpret this as $p=4$ for $q=\infty$)
\begin{equation}
\|r^{(n-1)/2}\varphi\|_{L_r^\infty L_\omega^p({\mathbb R}^n)}
\leq
\sqrt{p}
\|\partial_r \varphi\|_{L^2({\mathbb R}^n)}^{1/2}
\|\varphi\|_{L_r^2 L_\omega^q({\mathbb R}^n)}^{1/2}.
\end{equation}
Indeed, we obtain (2.7) directly from (2.19). 
Moreover, 
we obtain (2.11) immediately from 
the Sobolev embedding $W^{1,2}(S^2)\hookrightarrow L^q(S^2)$ 
for $2\leq q<\infty$ 
owing to the fact that in the condition for (2.20) to hold, 
the condition $2\leq p<4$ 
is equivalent to $2\leq q<\infty$. $\hfill\square$
\section{Weighted $L^2$-estimates}
\setcounter{equation}{0}
It is necessary to bound 
$M_4(u(t))$ by $N_4(u(t))$ 
for the completion of the energy integral argument 
(see Lemma 3.4 below). 
We carry out this by starting with the next crucial inequality 
due to Klainerman and Sideris $\cite{KS}$, the proof of which requires 
the use of the operator $S$; see $N_2(v(t))$ on the right-hand side of 
(3.1) below. 
In what follows we use the notation 
$\square_l:=\partial_t^2-c_l^2\Delta$. 
 
\smallskip

\noindent{\bf Lemma 3.1. (Klainerman--Sideris inequality)} 
{\it The inequality
\begin{equation}
M_2(v(t))
\leq
C
\bigl(
N_2(v(t))
+
\sum_{l=1}^N
t
\|\square_l v^l(t)\|_{L^2({\mathbb R}^n)}
\bigr)
\end{equation}
holds for any function $v=(v^1,\dots,v^N)$.
}

\smallskip

{\it Proof}. See Lemma 3.1 of $\cite{KS}$ 
(see also Lemma 7.1 of $\cite{ST}$). 
We have only to repeat essentially the same argument 
as in the proof of (3.1) of $\cite{KS}$. 
Note that the proof there is obviously valid for 
$n=2$ as well as $n=3$. $\hfill\square$

\smallskip

We also need the following auxiliary lemma, 
which compensates for the absence of 
$\partial_t^i$ $(i=2,3,4)$ in the norms 
appearing in (1.14), (2.5). 

\smallskip

\noindent{\bf Lemma 3.2.} {\it There exists a constant 
$\varepsilon^*>0$ depending on 
the propagation speeds $c_1,\dots,c_N$ 
and the coefficients 
on the right-hand side of $(1.7)$ or $(1.8)$ 
with the following property$:$ 
whenever a smooth solution $u=(u^1,\dots,u^N)$ 
to 
$(1.7)$ or $(1.8)$ satisfies 
\begin{equation}
\max\{\,|{\bar Z}^a\partial_\alpha u^i(t,x)|\,:\,
|a|\leq 1,\,0\leq\alpha\leq n,\,1\leq i\leq N\,\}
\leq
\varepsilon^*,
\end{equation}
the point-wise inequality
\begin{equation}
\sum_{i=1}^N
|\partial_t^2 u^i(t,x)|
\leq
C
\sum_{i=1}^N
\biggl(
\sum_{{1\leq m\leq n}\atop{0\leq\alpha\leq n}}
|\partial_{m\alpha}^2 u^i(t,x)|
+
\sum_{\alpha=0}^n
|\partial_\alpha u^i(t,x)|
\biggr)
\end{equation}
holds. 
Moreover, there holds for $|a|=1,2$
\begin{align}
&
\sum_{i=1}^N
|{\bar Z}^a\partial_t^2 u^i(t,x)|\\
&
\leq
C
\sum_{i=1}^N
\sum_{|b|=1}^{|a|}
\biggl(
\sum_{{1\leq m\leq n}\atop{0\leq\alpha\leq n}}
|{\bar Z}^b\partial_{m\alpha}^2 u^i(t,x)|
+
\sum_{\alpha=0}^n
|{\bar Z}^b\partial_\alpha u^i(t,x)|
\biggr).\nonumber
\end{align}
Here, $n=2,3$ for the solutions 
to $(1.8)$, $(1.7)$, respectively.
}

\smallskip

{\it Proof}. It suffices to prove the inequalities 
for the solutions to (1.7). The proof of the inequalities for 
the solutions to (1.8) is essentially the same. 
We first note the obvious equality for each $l=1,\dots,N$
\begin{align}
&
\partial_t^2 u^l
-
G_{ij}^{l,\alpha 00}
(\partial_\alpha u^i)\partial_t^2 u^j\\
&
=
c_l^2\Delta u^l
+
2G_{ij}^{l,\alpha m\gamma}(\partial_\alpha u^i)\partial_{m\gamma}^2 u^j
+
H_{ij}^{l,\alpha\beta}(\partial_\alpha u^i)\partial_\beta u^j.\nonumber
\end{align}
Whenever $|G_{ij}^{l,\alpha 00}\partial_\alpha u^i(t,x)|\leq 1/(2N)$ 
for any $j,l=1,\dots,N$, 
we get by (3.5)
\begin{align}
&
|\partial_t^2 u^l|
-
(2N)^{-1}
(|\partial_t^2 u^1|+\cdots+|\partial_t^2 u^N|)\\
&
\leq
c_l^2|\Delta u^l|
+
2|G_{ij}^{l,\alpha m\gamma}||\partial_\alpha u^i||\partial_{m\gamma}^2 u^j|
+
|H_{ij}^{l,\alpha\beta}||\partial_\alpha u^i||\partial_\beta u^j|.\nonumber
\end{align}
Summing both the sides of (3.6) over $l=1,\dots,N$, 
we see that the inequality (3.3) holds 
whenever $\max\{|\partial_\alpha u^i(t,x)|:\alpha=0,\dots,3,\,i=1,\dots,l\}$ 
is small enough. 

Next, let us prove (3.4). 
Using the commutation relations (2.1) and (2.4), we get for $|a|=1,2$
\begin{align}
&
\partial_t^2{\bar Z}^a u^l
-
G_{ij}^{l,\alpha 00}
(\partial_\alpha u^i)(\partial_t^2{\bar Z}^a u^j)\\
&
=
c_l^2\Delta{\bar Z}^a u^l
+
\sum_{{b+c=a}\atop{c\ne a}}
G_{ij}^{l,\alpha 00}
({\bar Z}^b\partial_\alpha u^i)
(\partial_t^2{\bar Z}^c u^j)\nonumber\\
&
+
\sum_{b+c=a}
2G_{ij}^{l,\alpha m\gamma}
({\bar Z}^b\partial_\alpha u^i)
({\bar Z}^c\partial_{m\gamma}^2 u^j)
+
\sum_{b+c=a}
H_{ij}^{l,\alpha\beta}
({\bar Z}^b\partial_\alpha u^i)
({\bar Z}^c\partial_\beta u^j).\nonumber
\end{align}
Noting the obvious fact 
$$
\sum_{{b+c=a}\atop{c\ne a}}
G_{ij}^{l,\alpha 00}
({\bar Z}^b\partial_\alpha u^i)
(\partial_t^2{\bar Z}^c u^j)
=
G_{ij}^{l,\alpha 00}
({\bar Z}^a\partial_\alpha u^i)
(\partial_t^2 u^j)
$$ 
for $|a|=1$, 
using (3.3) for the estimate of $|\partial_t^2 u^j(t,x)|$, 
and repeating the same argument as above, 
we see that the inequality (3.4) holds for $|a|=1$ 
whenever 
$\max\{\,|{\bar Z}^a\partial_\alpha u^i(t,x)|\,:\,
|a|\leq 1,\,0\leq\alpha\leq 3,\,1\leq i\leq N\,\}$ 
is small enough. 
Finally, using (3.3) and (3.4) with $|a|=1$ 
for the estimate of $|\partial_t^2{\bar Z}^c u^j(t,x)|$ 
$(|c|=0,1)$ (see the second term on the right-hand side of 
(3.7)) and repeating the above argument, 
we see that the inequality (3.4) holds for $|a|=2$ 
whenever 
$\max\{\,|{\bar Z}^a\partial_\alpha u^i(t,x)|\,:\,
|a|\leq 1,\,0\leq\alpha\leq 3,\,1\leq i\leq N\,\}$ 
is small enough. 
We have finished the proof.$\hfill\square$

\smallskip

Lemma 3.2 is useful in proving the following:

\noindent{\bf Lemma 3.3.} {\it 
Let $u=(u^1,\dots,u^N)$ be a smooth solution to 
$(1.7)$ or $(1.8)$ defined in $(0,T)\times{\mathbb R}^n$
satisfying
\begin{equation}
\sup_{(0,T)\times{\mathbb R}^n}
\max\{\,|{\bar Z}^a\partial_\alpha u^i(t,x)|\,:\,
|a|\leq 1,\,0\leq\alpha\leq n,\,1\leq i\leq N\,\}
\leq
\varepsilon^*.
\end{equation}
Then the following inequalities hold$:$ 
for each $l=1,2,\dots,N$ and $|a|\leq 2$,
\begin{equation}
t\|\square_l{\bar Z}^a u^l(t)\|_{L^2({\mathbb R}^3)}
\leq
CN_4^2(u(t))
+
CN_4(u(t))M_4(u(t)),\,\,0<t<T
\end{equation}
when $u$ is a solution to $(1.7)$, 
\begin{equation}
t\|
\square_l{\bar Z}^a u^l(t)
\|_{L^2({\mathbb R}^2)}
\leq 
CN_4^3(u(t))
+
CN_4^2(u(t))M_4(u(t)),\,\,0<t<T
\end{equation}
when $u$ is a solution to $(1.8)$.
}

\smallskip

{\it Proof}. We start with the $3$D case (3.9). 
Obviously, it suffices to deal with $|a|=2$. 
Taking account of the form of the quadratic nonlinear terms of (1.7) 
and using the commutation relations (2.1)--(2.4) and 
the point-wise inequality (3.4), we get
\begin{eqnarray}
& &
\|\square_l{\bar Z}^a u(t)\|_{L^2({\mathbb R}^3)}\\
& &
\leq
C
\sum_{i,j=1}^N
\sum_{{|b|+|c|}\atop{\leq 2}}
\bigl(
\|
(\partial{\bar Z}^b u^i(t))
\partial\partial_x{\bar Z}^c u^j(t)
\|_{L^2}
+
\|
(\partial{\bar Z}^b u^i(t))
\partial{\bar Z}^c u^j(t)
\|_{L^2}
\bigr).\nonumber
\end{eqnarray}
(Here, and in the following as well, 
we use the notation $\partial$ 
to mean any of the standard partial differential operators 
$\partial_a$ ($a=0,\dots,n$).) 
It is enough to handle only the case 
$|b|+|c|=2$. 
We use the notation 
$B_i:=\{x\in{\mathbb R}^3:|x|<(c_i/2)t+1\}$, 
with $B'_i$ being the complement of $B_i$. 
Using the triangle inequality, 
(2.13), (2.14), (2.11), and the Sobolev embedding 
$W^{1,2}(S^2)\hookrightarrow L^{\infty-}(S^2)$, 
we get for each $i,j=1,\dots,N$
\begin{align}
&
\sum_{|b|+|c|=2}
\|
(\partial{\bar Z}^b u^i(t))
\partial\partial_x{\bar Z}^c u^j(t)
\|_{L^2({\mathbb R}^3)}\\
&
\leq
C\sum_{|c|=2}
\langle t\rangle^{-1}
\bigl(
\|\langle c_it-r\rangle\partial u^i(t)\|_{L^\infty(B_i)}
\|\partial\partial_x{\bar Z}^c u^j(t)\|_{L^2({\mathbb R}^3)}\nonumber\\
&
\hspace{2.5cm}
+
\|r\partial u^i(t)\|_{L^\infty(B_i')}
\|\partial\partial_x{\bar Z}^c u^j(t)\|_{L^2({\mathbb R}^3)}
\bigr)\nonumber\\
&
+
C\sum_{|b|=|c|=1}
\langle t\rangle^{-1}
\bigl(
\|\langle c_it-r\rangle\partial{\bar Z}^b u^i(t)\|_{L^\infty(B_i)}
\|\partial\partial_x{\bar Z}^c u^j(t)\|_{L^2({\mathbb R}^3)}\nonumber\\
&
\hspace{2.5cm}
+
\|r\partial{\bar Z}^b u^i(t)\|_{L_r^\infty L_\omega^{2+}(B_i')}
\|\partial\partial_x{\bar Z}^c u^j(t)\|
_{L_r^2L_\omega^{\infty-}({\mathbb R}^3)}
\bigr)\nonumber\\
&
+
C\sum_{|b|=2}
\langle t\rangle^{-1}
\bigl(
\|\partial{\bar Z}^b u^i(t)\|_{L^2({\mathbb R}^3)}
\|\langle c_jt-r\rangle\partial\partial_x u^j(t)\|_{L^\infty(B_j)}\nonumber\\
&
\hspace{2.5cm}
+
\|\partial{\bar Z}^b u^i(t)\|_{L^2({\mathbb R}^3)}
\|r\partial\partial_x u^j(t)\|_{L^\infty(B_j')}
\bigr)\nonumber\\
&
\leq
C\langle t\rangle^{-1}
\bigl(
N_4(u(t))+M_4(u(t))
\bigr)
N_4(u(t)).\nonumber
\end{align}
(Here and in what follows, 
by $2+$ and $\infty-$ we mean 
arbitrary numbers $p_2$ and $p_3$, respectively, 
such that 
$p_2>2$, $p_3<\infty$, and $1/2=1/p_2+1/p_3$.) 
For the second norm on the right-hand side of (3.11), 
we obtain in the same way as above 
\begin{align}
&
\sum_{|b|+|c|=2}
\|
(\partial{\bar Z}^b u^i(t))
\partial{\bar Z}^c u^j(t)
\|_{L^2({\mathbb R}^3)}
\\
&
\leq
C\langle t\rangle^{-1}
\bigl(
N_4(u(t))
+
M_4(u(t))
\bigr)
N_4(u(t)).\nonumber
\end{align}
We have finished the proof of (3.9).

Next, let us prove (3.10). 
Again, we have only to handle the case $|a|=2$. 
As in (3.11), we get
\begin{align}
\|\square_l{\bar Z}^a u^l(t)\|_{L^2({\mathbb R}^2)}
\leq
C
\sum
\bigl(
&\|
(\partial{\bar Z}^b u^i(t))
(\partial{\bar Z}^c u^j(t))
\partial\partial_x{\bar Z}^d u^k(t)\|_{L^2}\\
&+
\|(\partial{\bar Z}^b u^i(t))
(\partial{\bar Z}^c u^j(t))
\partial{\bar Z}^d u^k(t)\|_{L^2}
\bigr).\nonumber
\end{align}
Here, the sum has been taken over 
$i,j,k=1,\dots,N$ and 
$b,c,d$ with 
$|b|+|c|+|d|\leq 2$. 
We must treat the two cases 
$|d|=0$ and $|d|=1,2$ separately. 
In the case $|d|=0$, 
assuming $|b|\leq 1$ and $|c|\leq 2$ 
without loss of generality, 
we obtain 
by (2.12), (2.14) 
\begin{align}
&
\|
(\partial{\bar Z}^b u^i(t))
(\partial{\bar Z}^c u^j(t))
\partial\partial_x u^k(t)
\|_{L^2({\mathbb R}^2)}\\
&
\leq
C\langle t\rangle^{-1}
\|
\langle c_i t-r\rangle^{1/2}
\partial{\bar Z}^b u^i(t)
\|_{L^\infty (B_{i,k})}\nonumber\\
&
\hspace{1cm}
\times
\|
\partial{\bar Z}^c u^j(t)
\|_{L^2({\mathbb R}^2)}
\|
\langle c_k t-r\rangle^{1/2}
\partial\partial_x u^k(t)
\|_{L^\infty(B_{i,k})}\nonumber\\
&
+
C\langle t\rangle^{-1}
\|
r^{1/2}
\partial{\bar Z}^b u^i(t)
\|_{L^\infty(B_{i,k}')}\nonumber\\
&
\hspace{1cm}
\times
\|
\partial{\bar Z}^c u^j(t)
\|_{L^2({\mathbb R}^2)}
\|
r^{1/2}
\partial\partial_x u^k(t)
\|_{L^\infty(B_{i,k}')}\nonumber\\
&
\leq
C\langle t\rangle^{-1}
\bigl(
N_4(u(t))+M_4(u(t))
\bigr)
N_4^2(u(t)),\nonumber
\end{align}
where 
$B_{i,k}:=\{x\in{\mathbb R}^2:|x|<\min\{c_it/2,c_kt/2\}+1\}$ 
with $B_{i,k}'$ being its complement. 
On the other hand, for $|d|=1,2$, 
we obtain by assuming 
$|b|=0$, $|c|\leq 1$ 
without loss of generality
\begin{align}
&
\|
(\partial u^i(t))
(\partial{\bar Z}^c u^j(t))
\partial\partial_x{\bar Z}^d u^k(t)
\|_{L^2({\mathbb R}^2)}\\
&
\leq
\|
(\partial u^i(t))
\partial{\bar Z}^c u^j(t)
\|_{L^\infty({\mathbb R}^2)}
\|
\partial\partial_x{\bar Z}^d u^k(t)
\|_{L^2({\mathbb R}^2)}\nonumber\\
&
\leq
C\langle t\rangle^{-1}
\bigl(
\|
\langle c_i t-r\rangle^{1/2}\partial u^i(t)
\|_{L^\infty(B_{i,j})}
\|
\langle c_j t-r\rangle^{1/2}\partial{\bar Z}^c u^j(t)
\|_{L^\infty(B_{i,j})}\nonumber\\
&
\hspace{1cm}
+
\|
r^{1/2}\partial u^i(t)
\|_{L^\infty(B_{i,j}')}
\|
r^{1/2}\partial{\bar Z}^c u^j(t)
\|_{L^\infty(B_{i,j}')}
\bigr)
\|
\partial\partial_x{\bar Z}^d u^k(t) 
\|_{L^2({\mathbb R}^2)}\nonumber\\
&
\leq
C\langle t\rangle^{-1}
\bigl(
N_4(u(t))+M_4(u(t))
\bigr)
N_4^2(u(t)).\nonumber
\end{align}
By (3.15)--(3.16), we have obtained the desired estimate 
of the first term on the right-hand side of (3.14). 
The proof of the estimate for the second term on its right-hand side 
is quite similar. We may omit it. 
The proof of Lemma 3.3 has been finished.$\hfill\square$

\smallskip

\noindent{\bf Lemma 3.4.} {\it 
There exists a small, positive constant 
$\delta_0$ with the following property$:$ 
suppose that, 
for a local smooth solution $u$ of $(1.7)$ or $(1.8)$, 
the supremum of $N_4(u(t))$ over an interval 
$(0,T)$ is sufficiently small so that 
\begin{equation}
\sup_{0<t<T}N_4(u(t))
\leq
\delta_0
\end{equation}
may hold. 
Then, the inequality
\begin{equation}
M_4(u(t))
\leq
CN_4(u(t)),\,\,0<t<T
\end{equation}
holds with a constant $C$ independent of $T$. 
}

\smallskip

{\it Proof}. Let us denote by $\delta_*$ the supremum of $N_4(u(t))$ over 
the interval $(0,T)$. 
By the Sobolev embedding, 
we see that (3.8) is satisfied 
when $\delta_*$ is sufficiently small. 
Then, we see that Lemma 3.1 with 
$v={\bar Z}^a u$ 
($|a|\leq 2$) and Lemma 3.3 
imply for $0<t<T$ 
\begin{align}
M_4(u(t))
&\leq
CN_4(u(t))
+
CN_4^\nu(u(t))
\bigl(
N_4(u(t))+M_4(u(t))
\bigr)\\
&\leq
C(1+\delta_*^\nu)N_4(u(t))
+
C\delta_*^\nu M_4(u(t)),\nonumber
\end{align}
($\nu=1,2$ for (1.7), (1.8), respectively) 
from which we easily verify the existence of 
the constant $\delta_0$, as claimed in the lemma.
$\hfill\square$

\noindent{\it Remark}. In the above proof, 
especially when absorbing $C\delta_*^\nu M_4(u(t))$ 
into the left-hand side of (3.19), 
we have used the fact that 
$M_4(u(t))$ is finite for $t\in (0,T)$. 
Indeed, using (3.1), (3.17), and the standard Sobolev embedding, 
we get
$M_4(u(t))\leq
CN_4(u(t))+CtN_4^{\nu+1}(u(t))<\infty$.  
\section{Estimate for $N_4(u(t))$}
\setcounter{equation}{0}
For the given smooth and compactly supported initial data, 
let us assume (1.12) for a sufficiently 
small $\varepsilon_0>0$ such that $2\varepsilon_0\leq\delta_0$ 
(see Lemma 3.4 for $\delta_0$). 
By the local existence theorem 
mentioned in Section 1, 
a unique smooth solution exists locally in time. 
Note that it is compactly supported 
for fixed times by the finite speed of propagation. 
Let $T^*$ be the supremum of the set 
of all $T>0$ such that 
this solution to (1.7) is defined in 
$(0,T)\times{\mathbb R}^3$ 
and satisfies 
$$
\sup_{0<t<T}N_4(u(t))<\infty.
$$
When considering (1.8), 
we define $T^*$ in the same way. 
When $T^*=\infty$, nothing remains to be done. 
We may therefore suppose $T^*<\infty$. 

Recall the notation 
$\square_l=\partial_t^2-c_l^2\Delta$. 
We set $E_4(u(t))=N_4^2(u(t))$ 
(see (1.14) for the definition of $N_4(u(t))$). 
For (1.8), 
setting 
$Z^a=\partial_1^{a_1}\partial_2^{a_2}\Omega_{12}^{a_3}S^{a_4}$ 
for $a=(a_1,\dots,a_4)$ 
and letting $a_*$ stand for 
any multi-index $a=(a_1,\dots,a_4)$ with $a_4\leq 1$, 
we have the energy equality by the standard argument
\begin{eqnarray}
& &
{\tilde E}'_4(u(t))\\
& &
=
\sum_{{1\leq l\leq N}\atop{|a_*|=3}}
\int_{{\Bbb R}^2}
G_{ijk}^{l,\alpha\beta\gamma\delta}
(\partial_\alpha u^i)
(\partial_\beta u^j)
([Z^{a_*},\partial_\gamma\partial_\delta]u^k)
\partial_t Z^{a_*}u^l dx\nonumber\\
& &
+
\sum_{{1\leq l\leq N}\atop{|a_*|=3}}
\sum_{{b+c+d=a_*}\atop{d\ne a_*}}
\int_{{\Bbb R}^2}
G_{ijk}^{l,\alpha\beta\gamma\delta}
(Z^b\partial_\alpha u^i)
(Z^c\partial_\beta u^j)
(Z^d\partial_\gamma\partial_\delta u^k)
\partial_t Z^{a_*}u^l dx\nonumber\\
& &
+
\sum_{{1\leq l\leq N}\atop{|a_*|\leq 2}}
\sum_{b+c+d=a_*}
\int_{{\Bbb R}^2}
G_{ijk}^{l,\alpha\beta\gamma\delta}
(Z^b\partial_\alpha u^i)
(Z^c\partial_\beta u^j)
(Z^d\partial_\gamma\partial_\delta u^k)
\partial_t Z^{a_*}u^l dx\nonumber\\
& &
-
\sum_{{1\leq l\leq N}\atop{|a_*|=3}}
\int_{{\Bbb R}^2}
G_{ijk}^{l,\alpha\beta p\delta}
\bigl(
\partial_p
(
(\partial_\alpha u^i)
\partial_\beta u^j
)
\bigr)
(\partial_\delta Z^{a_*}u^k)
\partial_t Z^{a_*}u^l dx\nonumber\\
& &
-
\sum_{{1\leq l\leq N}\atop{|a_*|=3}}
\frac12\int_{{\Bbb R}^2}
\biggl(
G_{ijk}^{l,\alpha\beta 00}
\bigl(
\partial_t
(
(\partial_\alpha u^i)
\partial_\beta u^j
)
\bigr)
(\partial_t Z^{a_*}u^k)
\partial_t Z^{a_*}u^l \nonumber\\
& &
\hspace{2cm}
-G_{ijk}^{l,\alpha\beta pq}
\bigl(
\partial_t
(
(\partial_\alpha u^i)
\partial_\beta u^j
)
\bigr)
(\partial_q Z^{a_*}u^k)
\partial_p Z^{a_*}u^l
\biggr)dx\nonumber\\
& &
+
\sum_{{1\leq l\leq N}\atop{|a_*|\leq 3}}
\sum_{b+c+d=a_*}
\int_{{\Bbb R}^2}
H_{ijk}^{l,\alpha\beta\gamma}
(Z^b\partial_\alpha u^i)
(Z^c\partial_\beta u^j)
(Z^d\partial_\gamma u^k)
\partial_t Z^{a_*}u^l dx\nonumber\\
& &
+
\sum_{{1\leq l\leq N}\atop{|a_*|\leq 3}}
\int_{{\Bbb R}^2}
([\square_l, Z^{a_*}]u^l)
\partial_t Z^{a_*}u^l
dx,\,\,0<t<T^*.\nonumber
\end{eqnarray}
Here, 
taking into account the quasi-linear character of (1.8), 
we have introduced the modified energy
\begin{align}
{\tilde E}_4(u(t))&:=E_4(u(t))\\
&
-
\sum_{{|a_*|=3}\atop{1\leq l\leq N}}
\frac12
\int_{{\Bbb R}^2}
\biggl(
G_{ijk}^{l,\alpha\beta 00}
(\partial_\alpha u^i)
(\partial_\beta u^j)
(\partial_t Z^{a_*}u^k)
\partial_t Z^{a_*}u^l\nonumber\\
&
\hspace{3cm}
-
G_{ijk}^{l,\alpha\beta pq}
(\partial_\alpha u^i)
(\partial_\beta u^j)
(\partial_q Z^{a_*}u^k)
\partial_p Z^{a_*}u^l
\biggr)dx.\nonumber
\end{align}
Note that, in (4.1)--(4.2), repeated indices have been summed 
when lowered and uppered. 
Precisely, the Greek indices $\alpha,\beta,\gamma$, and $\delta$ 
run from $0$ to $2$, while the roman $p$ and $q$ from $1$ to $2$. 
The roman indices $i$, $j$, and $k$ run from $1$ to $N$. 

Similarly, for (1.7), 
setting 
$Z^a
=
\partial_1^{a_1}\partial_2^{a_2}\partial_3^{a_3}
\Omega_{12}^{a_4}\Omega_{13}^{a_5}\Omega_{23}^{a_6}S^{a_7}$ 
for $a=(a_1,\dots,a_7)$ 
and 
letting $a^*$ stand for 
any multi-index $a=(a_1,\dots,a_7)$ 
with $a_7\leq 1$, we get 
\begin{eqnarray}
& &
{\tilde E}'_4(u(t))\\
& &
=
\sum_{1\leq l\leq N}
\sum_{|a^*|=3}
\int_{{\Bbb R}^3}
G_{ij}^{l,\alpha\beta\gamma}
(\partial_\alpha u^i)
([Z^{a^*},\partial_\beta\partial_\gamma]u^j)
\partial_t Z^{a^*}u^l dx\nonumber\\
& &
+
\sum_{1\leq l\leq N}
\sum_{|a^*|=3}
\sum_{{b+c=a^*}\atop{c\ne a^*}}
\int_{{\Bbb R}^3}
G_{ij}^{l,\alpha\beta\gamma}
(Z^b\partial_\alpha u^i)
(Z^c\partial_\beta\partial_\gamma u^j)
\partial_t Z^{a^*}u^l dx\nonumber\\
& &
+
\sum_{{1\leq l\leq N}\atop{|a^*|\leq 2}}
\sum_{b+c=a^*}
\int_{{\Bbb R}^3}
G_{ij}^{l,\alpha\beta\gamma}
(Z^b\partial_\alpha u^i)
(Z^c\partial_\beta\partial_\gamma u^j)
\partial_t Z^{a^*}u^l dx\nonumber\\
& &
-
\sum_{1\leq l\leq N}
\sum_{|a^*|=3}
\int_{{\Bbb R}^3}
G_{ij}^{l,\alpha p\gamma}
(
\partial_p\partial_\alpha u^i
)
(\partial_\gamma Z^{a^*}u^j)
\partial_t Z^{a^*}u^l dx\nonumber\\
& &
-
\sum_{1\leq l\leq N}
\sum_{|a^*|=3}
\int_{{\Bbb R}^3}
\frac12
\bigl(
G_{ij}^{l,\alpha 00}
(\partial_t\partial_\alpha u^i)
(\partial_t Z^{a^*}u^j)
\partial_t Z^{a^*}u^l\nonumber\\
& &
\hspace{4.7cm}
-
G_{ij}^{l,\alpha pq}
(\partial_t\partial_\alpha u^i)
(\partial_q Z^{a^*}u^j)
\partial_p Z^{a^*}u^l
\bigr)dx\nonumber\\
& &
+
\sum_{1\leq l\leq N}
\sum_{|a^*|\leq 3}
\sum_{b+c=a^*}
\int_{{\Bbb R}^3}
H_{ij}^{l,\alpha\beta}
(Z^b\partial_\alpha u^i)
(Z^c\partial_\beta u^j)
\partial_t Z^{a^*}u^l dx\nonumber\\
& &
+
\sum_{1\leq l\leq N}
\sum_{|a^*|\leq 3}
\int_{{\Bbb R}^3}
([\square_l, Z^{a^*}]u^l)
\partial_t Z^{a^*}u^l
dx,\,\,0<t<T^*.\nonumber
\end{eqnarray}
Here, we have defined 
\begin{eqnarray}
& &
{\tilde E}_4(u(t))\\
& &
:=E_4(u(t))
-
\sum_{1\leq l\leq N}
\sum_{|a^*|=3}
\frac12
\int_{{\Bbb R}^3}
\bigl(
G_{ij}^{l,\alpha 00}
(\partial_\alpha u^i)
(\partial_t Z^{a^*}u^j)
\partial_t Z^{a^*}u^l \nonumber\\
& &
\hspace{6cm}
-G_{ij}^{l,\alpha pq}
(\partial_\alpha u^i)
(\partial_q Z^{a^*}u^j)
\partial_p Z^{a^*}u^l
\bigr)dx.\nonumber
\end{eqnarray}
As in (4.1)--(4.2), repeated indices have been summed in (4.3)--(4.4) 
when lowered and uppered. 
The Greek indices $\alpha,\beta$, and $\gamma$ 
run from $0$ to $3$, 
while the roman $p$ and $q$ from $1$ to $3$. 
The roman indices $i$ and  $j$ run from $1$ to $N$. 

We may suppose without loss of generality that, 
for ${\tilde E}_4(u(t))$ defined in (4.2), (4.4), 
the inequality
\begin{equation}
\frac23
E_4(u(t))
\leq
{\tilde E}_4(u(t))
\leq
\frac32
E_4(u(t))
\end{equation}
holds by the Sobolev embedding, 
whenever $N_4(u(t))$ is small enough. 
We also note that, 
owing to the commutation relations (2.3)--(2.4), 
the commutators $[Z^{a_*},\partial_\beta\partial_\gamma]$ 
and $[Z^{a^*},\partial_\beta\partial_\gamma]$, 
which appear in the first term 
on the right-hand side of (4.1), (4.3), 
have the form
\begin{align}
&[Z^{a_*},\partial_\beta\partial_\gamma]
=
\sum_{|b_*|\leq 2}
\sum_{\alpha,\delta=0}^2
C^{a_*,b_*}_{\beta\gamma\alpha\delta}\partial_\alpha\partial_\delta Z^{b_*},\\
&[Z^{a^*},\partial_\beta\partial_\gamma]
=
\sum_{|b^*|\leq 2}
\sum_{\alpha,\delta=0}^3
C^{a^*,b^*}_{\beta\gamma\alpha\delta}\partial_\alpha\partial_\delta Z^{b^*},\nonumber
\end{align}
respectively, 
for each $a_*,a^*$ $(|a_*|=|a^*|=3)$, 
$\beta$, and $\gamma$. 
Here, by $C^{a_*,b_*}_{\beta\gamma\alpha\delta}$ 
and 
$C^{a^*,b^*}_{\beta\gamma\alpha\delta}$, 
we mean suitable constants depending on 
$a_*$, $b_*$, 
$a^*$, $b^*$, and 
$\alpha$, $\beta$, $\gamma$, and $\delta$. 
(Note that by $b_*$ and $b^*$, we mean 
any multi-index 
$(b_1,\dots,b_4)$ with $b_4\leq 1$, 
and $(b_1,\dots,b_7)$ with 
$b_7\leq 1$, respectively.) 
Note also that, thanks to (2.1) and $a_4\leq 1$, 
the commutator $[\square_l, Z^{a_*}]$ appearing 
in the last term on the right-hand side of (4.1) 
is $0$ or $2\square_l$. 
A similar note applies equally to $[\square_l, Z^{a^*}]$ 
which appears on the right-hand side of (4.3). 

Let us start the estimate of $E_4(u(t))$ with the case $n=2$. 
We remark that, under the assumption that 
$N_4(u(t))$ is small enough, 
we have 
by repeating quite the same argument as in the proof of Lemma 3.2
\begin{equation}
\sum_{i=1}^N
|\partial_t^2 Z^a u^i(t,x)|
\leq
C\sum
\bigl(
|\partial_{m\alpha}^2 Z^b u^i(t,x)|
+
|\partial_\alpha Z^b u^i(t,x)|
\bigr)
\end{equation}
for any multi-index $a$ with 
$|a|\leq 2$, $a_4\leq 1$. 
Here, on the right-hand side 
the sum is taken over all 
$1\leq i\leq N$, 
$1\leq m\leq 2$, 
$0\leq\alpha\leq 2$, 
and $b$ with 
$b_k\leq a_k$ 
for all $1\leq k\leq 4$. 

In what follows, 
we assume $N_4(u(t))$ is small so that 
(4.7) may hold. 
Using the energy equality (4.1) and the commutation relations (2.3), 
(2.4), and (4.6)--(4.7), we get
\begin{eqnarray}
& &
{\tilde E}'_4(u(t))
\leq
C
\sum
\|(\partial Z^b u^i)
(\partial Z^c u^j)
\partial\partial_x Z^d u^k\|_{L^2({\mathbb R}^2)}
N_4(u)\\
& &
\hspace{1.7cm}
+
C
\sum
\|(\partial Z^b u^i)
(\partial Z^c u^j)
\partial Z^d u^k\|_{L^2({\mathbb R}^2)}
N_4(u).\nonumber
\end{eqnarray}
Here, on the first term on the right-hand side above 
the sum is taken over all $i,j,k=1,\dots,N$, 
and $b,c,d$ with 
$|b|+|c|+|d|\leq 3$, $|d|\leq 2$, 
$b_4+c_4+d_4\leq 1$. 
On the second term, the sum is taken over all 
 $i,j,k=1,\dots,N$, 
and $b,c,d$ with 
$|b|+|c|+|d|\leq 3$, $b_4+c_4+d_4\leq 1$. 
Obviously, we have only to deal with the case 
$|b|+|c|+|d|=3$. 
Moreover, we may focus on 
the case of $b_4+c_4+d_4=1$ because 
the argument otherwise becomes much simpler. 
When treating the second term on the right-hand side above, 
we may also assume 
$|b|\leq |c|\leq |d|$ 
(hence $|b|, |c|\leq 1$) 
without loss of generality. 
We treat the two cases $d_4=0$ and $d_4=1$, 
separately. 
When $d_4=0$, we know $b_4+c_4=1$ and $|d|\leq 2$. 
If $b_4=0$, then we get by (2.8), (2.14)
\begin{align}
&
\|
(\partial Z^b u^i)
(\partial Z^c u^j)
\partial Z^d u^k
\|_{L^2}\\
&
\leq
C\langle t\rangle^{-1}
\|
\langle c_i t-r\rangle^{1/2}
\partial Z^b u^i
\|_{L^4(B_{i,k})}
\|\partial Z^c u^j\|_{L^\infty}
\|
\langle c_k t-r\rangle^{1/2}
\partial Z^d u^k
\|_{L^4(B_{i,k})}\nonumber\\
&
+
C\langle t\rangle^{-1}
\|r^{1/2}\partial Z^b u^i\|_{L^\infty(B_{i,k}')}
\|r^{1/2}\partial Z^c u^j\|_{L^\infty(B_{i,k}')}
\|\partial Z^d u^k\|_{L^2}\nonumber\\
&
\leq
C\langle t\rangle^{-1}
N_4^3(u(t)).\nonumber
\end{align}
If $b_4=1$, then 
we know $c_4=0$ and thus we obtain the same bound as (4.9) 
by considering 
$\|\partial Z^b u^i\|_{L^\infty}
\|\langle c_jt-r\rangle^{1/2}\partial Z^c u^j\|_{L^4(B_{j,k})}$ 
in place of 
$\|
\langle c_i t-r\rangle^{1/2}
\partial Z^b u^i
\|_{L^4(B_{i,k})}
\|\partial Z^c u^j\|_{L^\infty}
$. 

When $d_4=1$, we know $b_4=c_4=0$ and thus obtain 
by (2.12), (2.14)
\begin{align}
&
\|
(\partial Z^b u^i)
(\partial Z^c u^j)
\partial Z^d u^k
\|_{L^2}\\
&
\leq
C\langle t\rangle^{-1}
\|
\langle c_i t-r\rangle^{1/2}
\partial Z^b u^i
\|_{L^\infty(B_{i,j})}
\|
\langle c_j t-r\rangle^{1/2}
\partial Z^c u^j
\|_{L^\infty(B_{i,j})}
\|
\partial Z^d u^k
\|_{L^2}\nonumber\\
&
+
C\langle t\rangle^{-1}
\|r^{1/2}\partial Z^b u^i\|_{L^\infty(B_{i,j}')}
\|r^{1/2}\partial Z^c u^j\|_{L^\infty(B_{i,j}')}
\|\partial Z^d u^k\|_{L^2}\nonumber\\
&
\leq
C\langle t\rangle^{-1}
N_4^3(u(t)).\nonumber
\end{align}
Next, let us consider the bound for the first term 
on the right-hand side of (4.8). 
We may suppose $|b|\leq |c|$ (thus $|b|\leq 1$) 
without loss of generality. 
We discuss the two cases $d_4=0$ and $d_4=1$, 
separately. 
In the former case, 
we further treat the two cases 
$|d|\leq 1$ and $|d|=2$, separately. 

Suppose $d_4=0$ and $|d|\leq 1$. 
If $b_4=1$, then 
we have $c_4=0$, $|c|\leq 2$ and thus obtain 
by (2.7) and the Sobolev embedding on $S^1$
\begin{align}
&
\|
(\partial Z^b u^i)
(\partial Z^c u^j)
\partial\partial_x Z^d u^k
\|_{L^2}\\
&
\leq
C\langle t\rangle^{-1}
\|\partial Z^b u^i\|_{L^\infty}
\|
\langle c_jt-r\rangle^{1/2}
\partial Z^c u^j
\|_{L^4(B_{j,k})}
\|
\langle c_kt-r\rangle^{1/2}
\partial\partial_x Z^d u^k
\|_{L^4(B_{j,k})}\nonumber\\
&
+
C\langle t\rangle^{-1}
\|
r^{1/2}\partial Z^b u^i
\|_{L^\infty(B_{j,k}')}
\|
r^{1/2}\partial Z^c u^j
\|_{L_r^\infty L^2_\omega(B_{j,k}')}
\|
\partial\partial_x Z^d u^k
\|_{L^2_r L^\infty_\omega}\nonumber\\
&
\leq
C\langle t\rangle^{-1}
N_4^3(u(t)).\nonumber
\end{align}
If $b_4=0$, then $c_4=1$ and obtain
\begin{align}
&
\|
(\partial Z^b u^i)
(\partial Z^c u^j)
\partial^2 Z^d u^k
\|_{L^2}\\
&
\leq
C\langle t\rangle^{-1}
\|
\langle c_it-r\rangle^{1/2}
\partial u^i
\|_{L^\infty(B_{i,k})}
\|
\partial Z^c u^j
\|_{L^2}
\|
\langle c_kt-r\rangle^{1/2}
\partial\partial_x  u^k
\|_{L^\infty(B_{i,k})}\nonumber\\
&
+
C\langle t\rangle^{-1}
\|
r^{1/2}\partial u^i
\|_{L^\infty(B_{i,k}')}
\|
\partial Z^c u^j
\|_{L^2}
\|
r^{1/2}
\partial\partial_x u^k
\|_{L^\infty(B_{i,k}')}\nonumber\\
&
\leq
C\langle t\rangle^{-1}
N_4^3(u(t))\nonumber
\end{align}
for $|c|=3$ (thus $|b|=0$), and 
\begin{align}
&
\|
(\partial Z^b u^i)
(\partial Z^c u^j)
\partial^2 Z^d u^k
\|_{L^2}\\
&
\leq
C\langle t\rangle^{-1}
\|
\langle c_it-r\rangle^{1/2}
\partial Z^b u^i
\|_{L^\infty(B_{i,k})}
\|
\partial Z^c u^j
\|_{L^4}
\|
\langle c_kt-r\rangle^{1/2}
\partial\partial_x Z^d u^k
\|_{L^4(B_{i,k})}\nonumber\\
&
+
C\langle t\rangle^{-1}
\|
r^{1/2}\partial Z^b u^i
\|_{L^\infty(B_{i,k}')}
\|
r^{1/2}
\partial Z^c u^j
\|_{L_r^\infty L_\omega^2(B_{i,k}')}
\|
\partial\partial_x Z^d u^k
\|_{L_r^2 L_\omega^\infty}\nonumber\\
&
\leq
C\langle t\rangle^{-1}
N_4^3(u(t))\nonumber
\end{align}
for $|c|\leq 2$. 
(Recall that we are assuming $d_4=0$, $|d|\leq 1$.)

Next, suppose $d_4=0$ and $|d|=2$. 
Then, we know 
$|b|=0$, $|c|=1$ 
(because of $|b|\leq |c|$ and 
$|b|+|c|+|d|=3$), 
and easily obtain by (3.18)
\begin{align}
&
\|
(\partial Z^b u^i)
(\partial Z^c u^j)
\partial^2 Z^d u^k
\|_{L^2}\\
&
\leq
C\langle t\rangle^{-1}
\|
\partial u^i
\|_{L^\infty}
\|
\partial Z^c u^j
\|_{L^\infty}
\|
\langle c_kt-r\rangle
\partial\partial_x Z^d u^k
\|_{L^2(B_k)}\nonumber\\
&
+
C\langle t\rangle^{-1}
\|
r^{1/2}\partial u^i
\|_{L^\infty(B_{k}')}
\|
r^{1/2}
\partial Z^c u^j
\|_{L^\infty(B_{k}')}
\|
\partial\partial_x Z^d u^k
\|_{L^2}\nonumber\\
&
\leq
C\langle t\rangle^{-1}
N_4^3(u(t)).\nonumber
\end{align}
(The definition of $B_i$ is given in the proof of Lemma 3.3.) 

Turn our attention to the case of $d_4=1$. 
We know $b_4=c_4=0$ and $|c|\leq 2$. 
We discuss the two cases $|d|=1$ and $|d|=2$, separately. 
If $|d|=1$, then we get
\begin{align}
&
\|
(\partial Z^b u^i)
(\partial Z^c u^j)
\partial^2 Z^d u^k
\|_{L^2}\\
&
\leq
C\langle t\rangle^{-1}
\|
\langle c_it-r\rangle^{1/2}
\partial Z^b u^i
\|_{L^\infty(B_{i,j})}
\|
\langle c_jt-r\rangle^{1/2}
\partial Z^c u^j
\|_{L^4(B_{i,j})}
\|
\partial\partial_x Z^d u^k
\|_{L^4}\nonumber\\
&
+
C\langle t\rangle^{-1}
\|
r^{1/2}\partial Z^b u^i
\|_{L^\infty(B_{i,j}')}
\|
r^{1/2}
\partial Z^c u^j
\|_{L_r^\infty L_\omega^2(B_{i,j}')}
\|
\partial\partial_x Z^d u^k
\|_{L_r^2 L_\omega^\infty}\nonumber\\
&
\leq
C\langle t\rangle^{-1}
N_4^3(u(t)).\nonumber
\end{align}
If $|d|=2$, then we know 
$|b|=0$, $|c|=1$ and thus easily obtain
\begin{align}
&
\|
(\partial Z^b u^i)
(\partial Z^c u^j)
\partial^2 Z^d u^k
\|_{L^2}\\
&
\leq
C\langle t\rangle^{-1}
\|
\langle c_it-r\rangle^{1/2}
\partial u^i
\|_{L^\infty(B_{i,j})}
\|
\langle c_jt-r\rangle^{1/2}
\partial Z^c u^j
\|_{L^\infty(B_{i,j})}
\|
\partial\partial_x Z^d u^k
\|_{L^2}\nonumber\\
&
+
C\langle t\rangle^{-1}
\|
r^{1/2}\partial u^i
\|_{L^\infty(B_{i,j}')}
\|
r^{1/2}
\partial Z^c u^j
\|_{L^\infty(B_{i,j}')}
\|
\partial\partial_x Z^d u^k
\|_{L^2}\nonumber\\
&
\leq
C\langle t\rangle^{-1}
N_4^3(u(t)).\nonumber
\end{align}
We have finished bounding the right-hand side of (4.8). 
Taking account of the equivalence between 
$E_4(u(t))$ and ${\tilde E}_4(u(t))$ (see (4.5)), 
we get from (4.9)--(4.16)
\begin{equation}
{\tilde E}'_4(u(t))
\leq
C\langle t\rangle^{-1}
N_4^2(u(t)){\tilde E}_4(u(t))
\end{equation}
as far as $N_4(u(t))$ is small enough. 

We turn our attention to the case $n=3$. 
In the same way as we got (4.8), 
we obtain by (4.3)
\begin{align}
{\tilde E}'_4(u(t))
&
\leq
C
\sum
\|(\partial Z^b u^i)
\partial\partial_x Z^c u^j\|_{L^2}
N_4(u)\\
&
+
C
\sum
\|(\partial Z^b u^i)
\partial Z^c u^j\|_{L^2}
N_4(u),\nonumber
\end{align}
where in the right-hand side, 
the sum is taken over 
all $i,j=1,\dots,N$ and 
$b,c$ with $|b|+|c|\leq 3$ 
($|c|\leq 2$ for the first term), 
$b_7+c_7\leq 1$. 
As in the case of $n=2$, 
it suffices to treat the terms with $|b|+|c|=3$ 
and $b_7+c_7=1$.

Let us first treat the second term on 
the right-hand side above. 
We may suppose $|b|\leq |c|$ 
(thus $|b|\leq 1$) 
without loss of generality. 
When $c_7=0$, we know $b_7=1$, $|c|\leq 2$ and thus obtain 
by (2.10), (2.11) and the Sobolev embedding on $S^2$
\begin{align}
\|
(\partial Z^b u^i)\partial Z^c u^j
\|_{L^2}
&
\leq
C\langle t\rangle^{-1}
\|
\partial Z^b u^i
\|_{L^3}
\|
\langle c_j t-r\rangle
\partial Z^c u^j
\|_{L^6(B_j)}\\
&
+
C\langle t\rangle^{-1}
\|
r\partial Z^b u^i
\|_{L_r^\infty L_\omega^{2+}(B_j')}
\|
\partial Z^c u^j
\|_{L_r^2 L_\omega^{\infty-}}\nonumber\\
&
\leq
C\langle t\rangle^{-1}
N_4^2(u(t)).\nonumber
\end{align}
When $c_7=1$, we know $b_7=0$ and thus obtain by (2.13), (2.14)
\begin{align}
&
\|
(\partial Z^b u^i)\partial Z^c u^j
\|_{L^2}\\
&
\leq
C\langle t\rangle^{-1}
\bigl(
\|
\langle c_i t-r\rangle
\partial Z^b u^i
\|_{L^\infty(B_i)}
+
\|
r
\partial Z^b u^i
\|_{L^\infty(B_i')}
\bigr)
\|
\partial Z^c u^j
\|_{L^2}\nonumber\\
&
\leq
C\langle t\rangle^{-1}
N_4^2(u(t)).\nonumber
\end{align}
Next, let us turn to the estimate of 
the first term on the right-hand side of (4.18). 
Again, we discuss the two cases 
$c_7=0$ and $c_7=1$, separately. 

Suppose $c_7=0$. 
We handle the two cases $|c|\leq 1$ and 
$|c|=2$, separately. 
When $c_7=0$ and $|c|\leq 1$, 
we know $b_7=1$ and thus obtain 
\begin{align}
&
\|
(\partial Z^b u^i)
\partial\partial_x Z^c u^j
\|_{L^2}\\
&
\leq
C\langle t\rangle^{-1}
\|
\partial Z^b u^i
\|_{L^2}
\bigl(
\|
\langle c_j t-r\rangle
\partial\partial_x u^j
\|_{L^\infty(B_j)}
+
\|
r
\partial\partial_x u^j
\|_{L^\infty(B_j')}
\bigr)\nonumber\\
&
\leq
C\langle t\rangle^{-1}
N_4^2(u(t))\nonumber
\end{align}
for $|b|=3$ (thus $|c|=0$), and 
\begin{align}
&
\|
(\partial Z^b u^i)
\partial\partial_x Z^c u^j
\|_{L^2}\\
&
\leq
C\langle t\rangle^{-1}
\|
\partial Z^b u^i
\|_{L^3}
\|
\langle c_j t-r\rangle
\partial\partial_x Z^c u^j
\|_{L^6(B_j)}\nonumber\\
&
+
C\langle t\rangle^{-1}
\|
r
\partial Z^b u^i
\|_{L_r^\infty L_\omega^{2+}(B_j')}
\|
\partial\partial_x Z^c u^j
\|_{L_r^2 L_\omega^{\infty-}}\nonumber\\
&
\leq
C\langle t\rangle^{-1}
N_4^2(u(t))\nonumber
\end{align}
for $|b|\leq 2$. When $c_7=0$ and $|c|=2$, 
we know $b_7=1$, $|b|\leq 1$ and 
thus easily get by (3.18)
\begin{align}
&
\|
(\partial Z^b u^i)
\partial\partial_x Z^c u^j
\|_{L^2}\\
&
\leq
C\langle t\rangle^{-1}
\|
\partial Z^b u^i
\|_{L^\infty}
\|
\langle c_j t-r\rangle
\partial\partial_x Z^c u^j
\|_{L^2(B_j)}\nonumber\\
&
+
C\langle t\rangle^{-1}
\|
r
\partial Z^b u^i
\|_{L^\infty(B_j')}
\|
\partial\partial_x Z^c u^j
\|_{L^2}\nonumber\\
&
\leq
C\langle t\rangle^{-1}
N_4^2(u(t)).\nonumber
\end{align}

Finally, suppose $c_7=1$. 
We know $b_7=0$ and $|b|\leq 2$. 
Let us discuss the two cases 
$|c|=1$ and $|c|=2$, separately. 
If $|c|=1$, then 
\begin{align}
&
\|
(\partial Z^b u^i)
\partial\partial_x Z^c u^j
\|_{L^2}\\
&
\leq
C\langle t\rangle^{-1}
\|
\langle c_i t-r\rangle
\partial Z^b u^i
\|_{L^6(B_i)}
\|
\partial\partial_x S u^j
\|_{L^3}\nonumber\\
&
+
C\langle t\rangle^{-1}
\|
r
\partial Z^b u^i
\|_{L_r^\infty L_\omega^{2+}(B_i')}
\|
\partial\partial_x S u^j
\|_{L_r^2 L_\omega^{\infty-}}\nonumber\\
&
\leq
C\langle t\rangle^{-1}
N_4^2(u(t)).\nonumber
\end{align}
If $|c|=2$, then 
we know $|b|\leq 1$ (and $b_7=0$) 
and thus easily get
\begin{align}
&
\|
(\partial Z^b u^i)
\partial\partial_x Z^c u^j
\|_{L^2}\\
&
\leq
C\langle t\rangle^{-1}
\bigl(
\|
\langle c_i t-r\rangle
\partial Z^b u^i
\|_{L^\infty(B_i)}
+
\|
r
\partial Z^b u^i
\|_{L^\infty(B_i')}
\bigr)
\|
\partial\partial_x Z^c u^j
\|_{L^2}\nonumber\\
&
\leq
C\langle t\rangle^{-1}
N_4^2(u(t)).\nonumber
\end{align}
We have finished the required estimates 
of the two terms on the right-hand side of (4.18). 
Combining (4.19)--(4.25) and recalling (4.5), we get
\begin{equation}
{\tilde E}'_4(u(t))
\leq
C\langle t\rangle^{-1}N_4(u(t)){\tilde E}_4(u(t))
\end{equation}
as far as $N_4(u(t))$ is small enough. 

We are in a position to complete the proof of Theorem 1.1. 
We prove Theorem 1.1 for the solutions to (1.7), 
because the proof for those to (1.8) is similar. 
Let $T_*$ be the supremum of the set 
of all $T>0$ such that 
this solution to (1.7) is defined in 
$(0,T)\times{\mathbb R}^3$ 
and small so that 
\begin{equation}
\sup_{0<t<T}N_4(u(t))\leq 2\varepsilon,
\end{equation}
where $\varepsilon:=N_4(u(0))$. 
By definition, we know $T_*\leq T^*$. 

If we assume 
\begin{equation}
\varepsilon\log(1+T_*)<B
\end{equation}
for the constant $B$ defined via 
$\exp\{C_1B\}=7/6$ (see (4.29) below for the constant $C_1>0$), 
then we will get a contradiction. 
Indeed, we get by (4.26)
\begin{equation}
{\tilde E}'_{4}(u(t))
\leq
2\varepsilon C_1(1+t)^{-1}{\tilde E}_4(u(t)),\,\,
0<t<T_*
\end{equation}
for a suitable constant $C_1>0$. This together with (4.5), 
(4.27)--(4.28) 
immediately yields
\begin{equation}
N_4(u(t))\leq \frac74\varepsilon<2\varepsilon,\,\,
0<t<T_*.
\end{equation}
We note that, thanks to the fact that 
$u(t,x)$ is smooth and compactly supported 
for fixed times, 
we easily see that 
\begin{equation}
N_4(u(t))\in C([0,T^*)).
\end{equation} 
(It is this simple proof of (4.31) that 
needs the smoothness of the solution 
and the compactness of the support for fixed times.) 
Recall $T_*\leq T^*$ by definition. 
If $T_*<T^*$, then in view of (4.31) 
the bound (4.30) obviously 
contradicts the definition of $T_*$. 
We thus see $T_*=T^*$. 
Recall that 
the system (1.7) is invariant under the translation 
of the time variable, and that 
the length of the interval of existence 
of $C^\infty$-solutions to (1.7) with data 
$(\varphi,\psi)$ given at $t=t_0$ 
depends only on the $H^3$-norm of 
$(\nabla\varphi,\psi)$ 
but it does not on $t_0$. 
Thanks to the bound 
$N_4(u(t))\leq 7\varepsilon/4$ 
$(0<t<T^*)$, 
we can therefore extend this solution $u(t,x)$ to a larger strip, 
say, 
$(0,T^*+T')\times{\mathbb R}^3$ 
(for some $T'>0$) with 
$$
\sup_{0<t<T^*+T'}
N_4(u(t))<\infty
$$
by solving (1.7) subject to 
the compactly supported $C^\infty$-data 
$(u(T^*-\delta,x),\partial_t u(T^*-\delta,x))$ 
given at $t=T^*-\delta$. 
(Here, by $\delta>0$, 
we mean a sufficiently small positive number.) 
This, however, contradicts the definition of $T^*$. 
We thus see that (4.28) is false and there holds
\begin{equation}
\varepsilon\log(1+T_*)\geq B, 
\end{equation}
by which we have finished the proof of Theorem 1.1 
for the solutions to (1.7).
\section{Proof of Theorem 1.2}
\setcounter{equation}{0}
The proof of Theorem 1.2 requires 
only obvious modifications of that of Theorem 1.1. 
We may therefore leave the details to the reader.

\bigskip

\noindent{\em Acknowledgement.} 
The author is grateful to Professors 
Hideo Kubo and Hideo Takaoka for inviting 
him to the conference ``Harmonic Analysis and Nonlinear Partial Differential
Equations" held in July 2016 at 
the Research Institute for Mathematical Sciences (RIMS), Kyoto University. 
He also thanks Professor Chengbo Wang for 
information on low regularity well-posedness 
for quasi-linear wave equations, 
which was very helpful to him. 
His thanks go as well to the referee for reading the manuscript carefully 
and making a valuable comment on the continuation argument for 
local solutions in Section 4. 
He was supported in part by 
the Grant-in-Aid for Scientific Research (C) (No.\,15K04955), 
Japan Society for the Promotion of Science (JSPS).

\begin{flushleft}
Kunio Hidano\\
Department of Mathematics\\
Faculty of Education\\
Mie University\\
1577 Kurima-machiya-cho, Tsu\\
Mie 514-8507 Japan\\
e-mail address:hidano@@edu.mie-u.ac.jp
\end{flushleft}
\end{document}